\documentclass[12pt,reqno]{amsart}

 \usepackage[margin=1in]{geometry}

\usepackage{amsrefs}
\usepackage{enumerate}
\usepackage[hidelinks,colorlinks]{hyperref}
\usepackage{amsmath}
\usepackage{amssymb}
\usepackage{mathtools}
\usepackage{mathrsfs}
\usepackage{color}

\usepackage[normalem]{ulem}
\usepackage{comment}
\usepackage{yfonts}
\usepackage{xcolor}

\newtheorem{theorem}{Theorem}[section]
\newtheorem{prop}[theorem]{Proposition}
\newtheorem{lemma}[theorem]{Lemma}

\newtheorem{definition}[theorem]{Definition}

\newcommand{\ZZ}{\mathbb{Z}}
\newcommand{\CC}{\mathbb{C}}

\newsavebox\foobox
\newlength{\foodim}

\definecolor{purple}{RGB}{80, 0, 130}

\newcommand{\defeq}{\vcentcolon=}

\usepackage{graphicx}

\newcommand\blfootnote[1]{%
  \begingroup
  \renewcommand\thefootnote{}\footnote{#1}%
  \addtocounter{footnote}{-1}%
  \endgroup
}

\hypersetup{
    linkcolor=red,
    citecolor=blue,
 }

\allowdisplaybreaks

\begin{document}

\title[A closed formula for linear recurrences with constant coefficients]{A closed formula for linear recurrences with constant coefficients}
\author[Bruda, Fang, Gilman, Marquez, Miller, Prapashtica, Son, Waheed, Wang]{Glenn Bruda, Bruce Fang, Pico Gilman, Raul Marquez, Steven J. Miller, Beni Prapashtica, Daeyoung Son, Saad Waheed, Janine Wang}

\blfootnote{
\begin{center}
\textit{Email addresses:}\texttt{~glenn.bruda@ufl.edu} (Glenn Bruda),\texttt{~fangbaojun2002@gmail.com} (Bruce Fang), \texttt{~pico@mit.edu} (Pico Gilman),\texttt{~raul.marquez02@utrgv.edu} (Raul Marquez),\texttt{~sjm1@williams.edu} (Steven J. Miller),\texttt{~beniprap@gmail.com} (Beni Prapashtica),\texttt{~ds15@williams.edu} (Daeyoung Son),\texttt{~sw21@uchicago.edu} (Saad Waheed),\texttt{~jjw3@williams.edu} (Janine Wang)
\end{center}
}

\begin{abstract} 
Given a linear recurrence of the form $c_n=a_1c_{n-1}+\cdots+a_j c_{n-j}$, it is well-known that $c_n=\sum_{r}p_r(n)r^n$, where the sum is taken over the set of characteristic roots and each $p_r(n)$ is some polynomial. We give a closed formula for the coefficients of each polynomial $p_r(n)$ for any linear recurrence of this form.
\end{abstract}

\maketitle

\section{Introduction}

Let $\{a_n\}_{1\leq n\leq j}$ be a sequence of complex numbers and $\{c_n\}_{n\geq0}$ be the sequence defined by the recurrence $c_n=a_{1}c_{n-1}+\cdots+a_j c_{n-j}$ for all $n\geq j$, where the initial conditions $\{c_k\}_{0\leq k\leq j-1}$ are fixed. Such a sequence $\{c_n\}$ is called a \emph{linear recurrence with constant coefficients}. The \emph{characteristic polynomial} of $\{c_n\}$ is the polynomial $p(z)=z^j-a_1 z^{j-1}-a_2 z^{j-2}-\cdots-a_{j-1} z-a_j$. A \emph{characteristic root} of the recurrence satisfied by $\{c_n\}$ is defined to be a root of $p$.

A classical result (e.g., \cite[Theorem IV.9]{analytic_combo}) in the theory of recurrence relations states that $\{c_n\}$ admits a \emph{generalized Binet formula}: for all $n\geq0$, we have that $c_n=\sum_{r}p_r(n)r^n$, where the sum is taken over the set of characteristic roots and each $p_r(n)$ is a polynomial depending on the initial conditions. We construct closed formulas for the coefficients of each $p_r(n)$.

\begin{theorem}\label{main_theorem}
    Define the sequence $\{c_n\}\subseteq\CC$ by the recurrence relation
    \begin{align}\label{c_n_recurrence}
        c_n \ = \ a_1c_{n-1}+\cdots+a_jc_{n-j}
    \end{align}
    with $c_0,c_1,\dots,c_{j-1}\in\CC$. For each characteristic root $r$ and nonnegative integer $k\leq j-1$, let
    \begin{enumerate}[(i)]
        \item $m_r$ be the multiplicity of $r$,
        \item $s_{r,k,h}\defeq\sum_{\ell\geq j-k}{k+\ell\choose h}a_{\ell}r^{-k-\ell}$, and
        \item ${N}_r\defeq(-1)^{m_r-1}\left(\sum_{\ell\geq m_r}{\ell\choose m_r}a_{\ell}r^{-\ell}\right)^{-1}$. 
    \end{enumerate}
    Then $c_n =\sum_r\sum_{i=1}^{m_r}\kappa_{r,i}n^{m_r-i}r^n$, where the outer sum is taken over the set of characteristic roots and
    \begin{align}\label{subset_formula_kappa}
        \kappa_{r,i} \ = \ \frac{(-1)^{i-1}}{(m_r-i)!}\sum_{\substack{\varnothing\neq\{d_1,\dots,d_p\}\subseteq\{0,1,\dots,i-1\}\\  d_1<\cdots<d_p=i-1}}{N}_r^p\sum_{k_1=0}^{j-1}c_{k_1}s_{r,k_1,d_1}\prod_{u=1}^{p-1}\sum_{k_2=0}^{j-1}\frac{k_2^{m_r-1-d_u}r^{k_2} s_{r,k_2,d_{u+1}}}{(-1)^{d_u-1}(m_r-1-d_u)!}.
    \end{align}
\end{theorem}

Computing $p_r(n)$ is traditionally done by algorithmically computing the generating function for $\{c_n\}$ and determining its partial fraction decomposition. Algorithms for computing partial fraction decompositions are well-studied, some of which can be found in \cite{horowitz,kung_and_tong,pfd_xin}.

From Theorem \ref{main_theorem}, one may obtain a closed formula for the partial fraction decomposition of a proper rational function, circumventing any algorithmic computations. A closed formula for such a partial fraction decomposition was also recently given by Chargeishvili, Fek\'esh\'azy, Somogyi, and Van Thurenhout in \cite{linapart}. Of course, from \cite[Equation 11]{linapart}, a closed formula for the linear recurrence may be obtained by applying Newton's binomial expansion to each term in the composition. However, our formula contrasts favorably with this approach in two ways: 
\begin{enumerate}[(i)]
    \item the formula given for each $\kappa_{r,i}$ is not a function of any of the other characteristic roots, and
    \item our formula provides an exact formula for each coefficient of $p_r(n)$ rather than representing $p_r(n)$ as a sum of binomial coefficients.
\end{enumerate}
Thus, if one wished to compute a particular $p_r(n)$, none of the roots other than $r$ would need to be computed. 

The primary application of Theorem \ref{main_theorem} is toward an implementation to symbolically compute the generalized Binet formula of any linear recurrence with constant coefficients. Most of the difficulty lies in accounting for multiple roots, since each leading coefficient $\kappa_{r,1}$ reduces to
\begin{align}
    \kappa_{r,1} \ = \ \frac{N_r}{(m_r-1)!}\sum_{k=0}^{j-1}c_k s_{r,k,0}.
\end{align}
Thus, Theorem \ref{main_theorem} simplifies dramatically if all of the characteristic roots of the recurrence are of multiplicity one. A prototypical class of recurrences with this property are the $j$-nacci numbers, which are defined by the recurrence $F(n,j)=\sum_{\ell=1}^{j}F(n-\ell,j)$ and initial conditions $F(0,j)=F(1,j)=\cdots=F(j-2,j)=0$ and $F(j-1,j)=1$. Letting $j\geq2$ and $n\geq0$, from Theorem \ref{main_theorem} we have the formula
    \begin{align}
        F(n,j) \ = \ \sum_r\frac{r^{n+1}(1-r^{-j})}{r^{j+1}-j},
    \end{align}
where the sum is taken over the set of roots of the characteristic polynomial $x^j-x^{j-1}-\cdots-x-1$. This follows from Theorem \ref{main_theorem} by a simplification of geometric series and the identity $\sum_{k=0}^{j-1}r^k=r^j$. This result is, of course, equivalent to the various closed forms for the $j$-nacci numbers given in the literature (e.g., \cite[Lemma 3.2]{closed_form_j_nacci}).

However, our formula permitting multiple roots allows us to consider more complicated recurrences. Exposited by Ahlgren in \cite{ahlgren}, Ramanujan stated in his last notebook \cite[pg.54]{last_notebook} that the circular summation of the sixth powers of the symmetric theta function decomposes into a product of Ramanujan's general theta function and another function $F_6$. By a formula for $F_6(q)$ given in \cite[Theorem 1]{ahlgren}, Zabolotskiy notes in \cite{ahlgren_sequence} that the exponents in the expansion of $F_6(q^2)$ are precisely the nonnegative integers congruent to $0,5,8,9$ $\mathrm{mod~} 12$, which satisfy the recurrence $c_n=2c_{n-1}-2c_{n-2}+2c_{n-3}-c_{n-4}$. By definition, $c_0=0$, $c_1=5$, $c_2=8$, and $c_3=9$. The characteristic roots are $1,i,-i$, which have multiplicities $2,1,1$ respectively. Applying Theorem \ref{main_theorem}, we find that $\kappa_{1,1}=3,~\kappa_{1,2}=1,~\kappa_{i,1}=(-1-i)/2,~\kappa_{-i,1}=(-1+i)/2$, thus obtaining
\begin{align}
    c_n \ = \ 3n+1+\frac{1}{2}(-1-i)i^n+\frac{1}{2}(-1+i)(-i)^n \ = \ 3n+1+\sin\left(\frac{\pi n}{2}\right)-\cos\left(\frac{\pi n}{2}\right).
\end{align}

We now consider another example demonstrating the power of allowing for multiple roots: the \emph{generalized pentagonal numbers} are the exponents in the expansion given by the celebrated pentagonal number theorem \cite{pentagonal_number_theorem}, which asserts the $q$-series identity
\begin{align}
    \prod_{n\geq1}(1-q^n) \ = \ \sum_{k\in\ZZ}(-1)^k q^{k(3k-1)/2}.
\end{align}
That is, $n$ is a generalized pentagonal number if and only if $n=k(3k-1)/2$ for some $k\in\ZZ$. We aim to prove a closed formula for the generalized pentagonal numbers in increasing order, which in turn would provide a singly infinite series expansion of monomials. By \cite{pentagonal_OEIS}, this sequence satisfies the recurrence $c_n=c_{n-1}+2c_{n-2}-2c_{n-3}-c_{n-4}+c_{n-5}$ with $c_0=0$, $c_1=1$, $c_2=2$, $c_3=5$, and $c_4=7$. The characteristic roots are $1$ and $-1$, which have multiplicities $3$ and $2$ respectively. Applying Theorem \ref{main_theorem}, we find that $\kappa_{1,1}=3/8,~\kappa_{1,2}=3/8,~\kappa_{1,3}=1/16,~\kappa_{-1,1}=-1/8,~\kappa_{-1,2}=-1/16$, thus obtaining
\begin{align}
    c_n \ = \ \frac{3n^2}{8}+\frac{3n}{8}+\frac{1}{16}+\left(-\frac{n}{8}-\frac{1}{16}\right)(-1)^n,
\end{align}
which matches the formula given by Barry in \cite[Section 7]{barry_formula}.

The rest of the paper is dedicated to proving Theorem \ref{main_theorem}, followed by a discussion on directions for future work.

\section{Proof of Theorem \ref{main_theorem}}

To prove Theorem \ref{main_theorem}, we appeal to the fact that the coefficient of the asymptotically dominant term in the expansion formula for $\{c_n\}$ is straightforward to compute by \cite[Note IV.26]{analytic_combo}. In particular, to compute the coefficient of $n^{m_r-i} r^n$, we subtract all terms which asymptotically dominate $n^{m_r-i} r^n$, obtaining a new linear recurrence $\{c_n'\}$ satisfying $c'_n=\kappa_{r,i}n^{m_r-i} r^n(1+O(1/n))$. To leverage \cite[Note IV.26]{analytic_combo}, we first must compute the generating function of $\{c_n\}$ in closed form. To do this, we break $\{c_n\}$ into a linear combination of a simpler sequence, which we call the $(a_1,\dots,a_j)$-nacci sequence.

\begin{definition}[$(a_1,\dots,a_j)$-nacci sequence]\label{S_nacci}
    Let $\{a_n\}_{1\leq n\leq j}$ be a sequence of complex numbers. Define the sequence $\{F_{n}^{(a_1,\dots,a_j)}\}_{n\geq0}$, which we call the $(a_1,\dots,a_j)$-nacci sequence, by the recurrence relation
    \begin{align}
        F_{n}^{(a_1,\dots,a_j)} \ = \ \sum_{\ell\geq 1} a_{\ell}F_{n-\ell}^{(a_1,\dots,a_j)},
    \end{align}
    with $F_{0}^{(a_1,\dots,a_j)}=F_{1}^{(a_1,\dots,a_j)}=\cdots=F_{j-2}^{(a_1,\dots,a_j)}=0$ and $F_{j-1}^{(a_1,\dots,a_j)}=1$.
\end{definition}

The upside of considering the simpler $(a_1,\dots,a_j)$-nacci sequence is that its generating function is easy to compute in closed form. Thus, once we precisely describe the linear combination yielding $\{c_n\}$, the desired generating function formula for $\{c_n\}$ quickly follows.

\begin{lemma}\label{thm:lem1}
    Let $\{a_n\}_{1\leq n\leq j}$ be a sequence of complex numbers. Then
    \begin{align}
        \sum_{n\geq0}F_n^{(a_1,\dots,a_j)}z^k \ = \ \frac{z^{j-1}}{1-\sum_{\ell\geq 1}a_{\ell}z^\ell}.
    \end{align}
\end{lemma}

\begin{proof}
To determine the ordinary generating function of $\{F^{(a_1,\dots,a_j)}_{n}\}$, we use the standard ansatz (see ~\cite{stackexchangeexample} for example) that this function is rational with denominator $1-\sum_{\ell\geq1}a_{\ell}z^\ell$. Indeed, we have
    \begin{align}
        \left(1-\sum_{\ell\geq1}a_{\ell}z^\ell\right)\sum_{k\geq0}F_k^{(a_1,\dots,a_j)}z^k \ &= \  \sum_{k\geq0}F_k^{(a_1,\dots,a_j)}z^k-\sum_{\ell\geq   1}\sum_{k\geq0}a_{\ell}F_k^{(a_1,\dots,a_j)}z^{k+\ell}\nonumber\\
        &= \ \sum_{k\geq0}F_k^{(a_1,\dots,a_j)}z^k-\sum_{\ell\geq1}\sum_{k\geq \ell}a_{\ell}F_{k-\ell}^{(a_1,\dots,a_j)}z^{k}\nonumber\\
        &\overset{*}{=} \ F_{j-1}^{(a_1,\dots,a_j)}z^{j-1}+\sum_{k\geq j}\left(F_{k}^{(a_1,\dots,a_j)}-\sum_{\ell\geq1}a_{\ell}F_{k-\ell}^{(a_1,\dots,a_j)}\right)z^k\nonumber \\ 
        &= \ F_{j-1}^{(a_1,\dots,a_j)}z^{j-1} \ = \ z^{j-1},
    \end{align}
    where the starred equality follows since $F_{k}^{(a_1,\dots,a_j)}=0$ if $k\leq j-2$. So
    \begin{align}
        \sum_{k\geq0}F_k^{(a_1,\dots,a_j)}z^k \ = \ \frac{z^{j-1}}{1-\sum_{\ell\geq1}a_{\ell}z^\ell},
    \end{align}
    as desired.
\end{proof}

\begin{lemma}\label{thm:mtp}
     Define the sequence $\{c_n\}\subseteq\CC$ by the recurrence relation $c_n = a_1c_{n-1}+a_2c_{n-2}+\cdots+a_jc_{n-j}$ with $c_0,c_1,\dots,c_{j-1}\in\CC$. Then for all $n\geq j$,
        \begin{align}
            c_n \ = \ \sum_{k=0}^{j-1} c_{k}{\sum_{\ell\geq j-k}a_{\ell}F_{n+j-1-k-\ell}^{(a_1,\dots,a_j)}}.
        \end{align}
\end{lemma}

In \cite[pg.6]{multiplicity_thm}, Liu proves a result similar to Lemma \ref{thm:mtp}, using instead an explicit formula for $F_{n}^{(a_1,\dots,a_j)}$. As we are ultimately only concerned with the generating function of $F_{n}^{(a_1,\dots,a_j)}$, we need not consider its explicit form. 

We now recall the definition of the Kronecker delta for use in the proof of Lemma \ref{thm:mtp}. For two numbers $k,k'$, the Kronecker delta is defined by
\begin{align}
    \delta_{k,k'} \ \defeq \
    \begin{cases}
        1 & \text{if~}k=k',\\
        0 & \text{if~}k\neq k'.
    \end{cases}
\end{align}

\begin{proof}[Proof of Lemma \ref{thm:mtp}]
    We have that
    \begin{align}
        Q(a_1,\dots,a_j) \ \defeq \ \begin{bmatrix} 
    a_1 & a_2 & \cdots & a_{j-1} & a_j\\
    1 & 0 & \cdots & 0 & 0 \\
    0 & 1 & \cdots & 0 & 0 \\
    \vdots & \vdots & \ddots & \vdots & \vdots \\
    0 & 0 & \cdots & 1 & 0
    \end{bmatrix}
    \end{align}
    is the $j\times j$ companion matrix of the $(a_1,\dots,a_j)$-nacci sequence. This matrix encodes the recurrence relation for the $(a_1,\dots,a_j)$-nacci sequence in that
\begin{align}
\begin{bmatrix}
    F^{(a_1,\dots,a_j)}_{n+1} \\
    F^{(a_1,\dots,a_j)}_n \\
    \vdots \\
    F^{(a_1,\dots,a_j)}_{n-j+2}
\end{bmatrix} \ = \ Q(a_1,\dots,a_j)\begin{bmatrix}
    F^{(a_1,\dots,a_j)}_n \\
    F^{(a_1,\dots,a_j)}_{n-1} \\
    \vdots \\
    F^{(a_1,\dots,a_j)}_{n-j+1}
\end{bmatrix}.
\end{align}
For each $0\leq k \leq j-1$, let $e_{k,n}$ be such that $c_n = \sum_{k=0}^{j-1} c_{k} e_{k,n}$. From the linear recurrence satisfied by $\{c_n\}$, we deduce $e_{k,n} = \sum_{\ell\geq1}a_{\ell}e_{k,n-\ell}$. As this is the same recurrence satisfied by the $(a_1,\dots,a_j)$-nacci sequence,
\begin{align}
\begin{bmatrix}
    e_{k,n+1} \\
    e_{k,n} \\
    \vdots \\
    e_{k,n-j+2}
\end{bmatrix} \ = \ Q(a_1,\dots,a_j)\begin{bmatrix}
    e_{k,n} \\
    e_{k,n-1} \\
    \vdots \\
    e_{k,n-j+1}
\end{bmatrix}.
\end{align}
Thus,

\begin{align}
    \begin{bmatrix}
    e_{k,n+j-1} \\
    e_{k,n+j-2} \\
    \vdots \\
    e_{k,n}
\end{bmatrix} \ = \ \left(Q(a_1,\dots,a_j)\right)^n\begin{bmatrix}
    e_{k,j-1} \\
    e_{k,j-2} \\
    \vdots \\
    e_{k,0}
\end{bmatrix} 
\ = \ \left(Q(a_1,\dots,a_j)\right)^n
\begin{bmatrix}
    \delta_{k,j-1} \\
    \delta_{k,j-2} \\
    \vdots \\
    \delta_{k,0}
\end{bmatrix}
.\label{lastentry}
\end{align}
We claim that $\left(Q(a_1,\dots,a_j)\right)^n$ equals
\begin{align*}
\begin{bmatrix}
   \sum_{\ell\geq1}a_{\ell}F_{n+j-1-\ell}^{(a_1,\dots,a_j)} & \sum_{\ell\geq2}a_{\ell}F_{n+j-\ell}^{(a_1,\dots,a_j)} & \cdots & \sum_{\ell\geq j-1}a_{\ell}F_{n+2j-3-\ell}^{(a_1,\dots,a_j)} & \sum_{\ell\geq j}a_{\ell}F_{n+2j-2-\ell}^{(a_1,\dots,a_j)}\\
    \vspace{1pt} & & & & \\
    \sum_{\ell\geq1}a_{\ell}F_{n+j-2-\ell}^{(a_1,\dots,a_j)} & \sum_{\ell\geq2}a_{\ell}F_{n+j-1-\ell}^{(a_1,\dots,a_j)} & \cdots & \sum_{\ell\geq j-1}a_{\ell}F_{n+2j-4-\ell}^{(a_1,\dots,a_j)} & \sum_{\ell\geq j}a_{\ell}F_{n+2j-3-\ell}^{(a_1,\dots,a_j)}\\
    \vdots & \vdots & \ddots & \vdots & \vdots \\
     \sum_{\ell\geq1}a_{\ell}F_{n+1-\ell}^{(a_1,\dots,a_j)} & \sum_{\ell\geq2}a_{\ell}F_{n+2-\ell}^{(a_1,\dots,a_j)} & \cdots & \sum_{\ell\geq j-1}a_{\ell}F_{n+j-1-\ell}^{(a_1,\dots,a_j)} & \sum_{\ell\geq j}a_{\ell}F_{n+j-\ell}^{(a_1,\dots,a_j)}\\
    \vspace{1pt} & & & & \\
     \sum_{\ell\geq1}a_{\ell}F_{n-\ell}^{(a_1,\dots,a_j)} & \sum_{\ell\geq2}a_{\ell}F_{n+1-\ell}^{(a_1,\dots,a_j)} & \cdots & \sum_{\ell\geq j-1}a_{\ell}F_{n+j-2-\ell}^{(a_1,\dots,a_j)} & \sum_{\ell\geq j}a_{\ell}F_{n+j-1-\ell}^{(a_1,\dots,a_j)}
\end{bmatrix}
\end{align*}
for any integer $n\geq j$. Indeed, let $q_{a,b}^{(n)}$ denote the $(a,b)$\textsuperscript{th} entry of $(Q(a_1,\dots,a_j))^n$. By ~\cite[Theorem 3.2]{powermatrix},
\begin{align}
    \sum_{n\geq0}q_{a,b}^{(n)}z^n \ = \ \frac{z^{a-b}\left(1-\sum_{\ell=1}^{b-1}a_{\ell}z^{\ell}\right)}{1-\sum_{\ell\geq1}a_{\ell}z^{\ell}}
\end{align}
if $a\geq b$, and
\begin{align}
    \sum_{n\geq0}q_{a,b}^{(n)}z^n \ = \ \frac{\sum_{\ell\geq b}a_{\ell}z^{a-b+\ell}}{1-\sum_{\ell\geq1}a_{\ell}z^{\ell}}
\end{align}
if $a<b$. Lemma \ref{thm:lem1} gives that
\begin{align}
    \sum_{n\geq0}F_n^{(a_1,\dots,a_j)}z^k \ = \ \frac{z^{j-1}}{1-\sum_{\ell\geq 1}a_{\ell}z^\ell}.
\end{align}
From these formulas, we see that
\begin{align}
    q_{a,b}^{(n)} \ = \ F_{n+j-1+b-a}^{(a_1,\dots,a_j)}-\sum_{\ell=1}^{b-1}a_{\ell}F_{n+j-1+b-a-\ell}^{(a_1,\dots,a_j)} \ = \ \sum_{\ell\geq b}F_{n+j-1+b-a-\ell}^{(a_1,\dots,a_j)}
\end{align}
for any $a,b$, verifying the formula for $(Q(a_1,\dots,a_j))^n$.

The last entry $e_{k,n}$, of the vector in \eqref{lastentry} is obtained by taking the $(j-k)$\textsuperscript{th} entry (with entries numbered from left to right starting from 1) of the last row of  $\left(Q(a_1,\dots,a_j)\right)^n$. Thus,
\begin{align}
e_{k,n} \ = \ \sum_{\ell\geq j-k}a_{\ell}F_{n+j-1-k-\ell}^{(a_1,\dots,a_j)},
\end{align}
completing the proof.
\end{proof}

From Lemmas \ref{thm:lem1} and \ref{thm:mtp}, we obtain a closed formula for the generating function of an arbitrary linear recurrence with constant coefficients. We note that this formula is equivalent to other such generating function formulas in the literature, such as \cite[Theorem 5]{gen_func_formula}; however, writing the formula in this particular manner is conducive to ultimately simplifying the formula for $\kappa_{r,i}$.

\begin{prop}\label{generating_function_formula}
    Define the sequence $\{c_n\}\subseteq\CC$ by the recurrence relation
    \begin{align}
        c_n \ = \ a_1c_{n-1}+a_2c_{n-2}+\cdots+a_jc_{n-j}
    \end{align}
    with $c_0,c_1,\dots,c_{j-1}\in\CC$. Then
    \begin{align}
        \sum_{n\geq0}c_n z^n \ = \ \frac{\sum_{k=0}^{j-1}c_k\sum_{\ell\geq j-k} a_{\ell}z^{k+\ell}+(1-a_1z-\cdots-a_jz^j)\sum_{k=0}^{j-1}c_k z^k}{1-a_1z-\cdots-a_jz^j}.
    \end{align}
\end{prop}

\begin{proof}
    By Lemmas \ref{thm:lem1} and \ref{thm:mtp},
    \begin{align}
        \sum_{n\geq0}c_n z^n \ &= \ \sum_{k=0}^{j-1}c_kz^k+\sum_{n\geq j}z^n\sum_{k=0}^{j-1} c_{k}{\sum_{\ell\geq j-k}a_{\ell}F_{n+j-1-k-\ell}^{(a_1,\dots,a_j)}}\nonumber\\
        &= \ \sum_{k=0}^{j-1}c_kz^k+\sum_{k=0}^{j-1} c_{k}{\sum_{\ell\geq j-k}a_{\ell}\sum_{n\geq j}F_{n+j-1-k-\ell}^{(a_1,\dots,a_j)}}z^n\nonumber\\
        &= \ \sum_{k=0}^{j-1}c_kz^k+\sum_{k=0}^{j-1} c_{k}{\sum_{\ell\geq j-k}a_{\ell}\sum_{n\geq 2j-1-k-\ell}F_{n}^{(a_1,\dots,a_j)}}z^{n-j+1+k+\ell}\nonumber\\
        &= \ \sum_{k=0}^{j-1}c_kz^k+\sum_{k=0}^{j-1} c_{k}{\sum_{\ell\geq j-k}a_{\ell}\sum_{n\geq j-1}F_{n}^{(a_1,\dots,a_j)}}z^{n-j+1+k+\ell}\nonumber\\
        &= \ \sum_{k=0}^{j-1}c_kz^k+\sum_{k=0}^{j-1} c_{k}\sum_{\ell\geq j-k}a_{\ell}\frac{z^{k+\ell}}{1-a_1z-\cdots-a_jz^j}\nonumber\\
        &= \ \frac{\sum_{k=0}^{j-1}c_k\sum_{\ell\geq j-k} a_{\ell}z^{k+\ell}+(1-a_1z-\cdots-a_jz^j)\sum_{k=0}^{j-1}c_k z^k}{1-a_1z-\cdots-a_jz^j}.
    \end{align}

\end{proof}

The final lemma before the proof of Theorem \ref{main_theorem}, Lemma \ref{zeroing_lemma}, is crucial toward simplifying our formula for $\kappa_{r,i}$. In particular, it is due to this result that our formula for $\kappa_{r,i}$ does not depend on any characteristic root other than $r$.

\begin{lemma}\label{zeroing_lemma}
    Let $p(z)=z^j-a_1 z^{j-1}-\cdots-a_j$. Then for any nonzero distinct roots $r_1,r_2$ of $p$,
    \begin{align}
        \sum_{k=0}^{j-1}k^{M_1} r_1^k\sum_{\ell\geq j-k}{k+\ell\choose M_2} a_{\ell} r_2^{-k-\ell} \ = \ 0
    \end{align}
    for any nonnegative integers $M_1\leq m_{r_1}-1$ and $M_2\leq m_{r_2}-1$, where $m_{r_1}$ and $m_{r_2}$ are the multiplicities of $r_1$ and $r_2$ respectively. 
\end{lemma}

\begin{proof}
    Rewriting the sum as
    \begin{align}
        \sum_{k=0}^{j-1}k^{M_1} r_1^k\sum_{\ell\geq j-k}  {k+\ell\choose M_2}a_{\ell} r_2^{-k-\ell} \ &= \ \sum_{\ell=1}^{j}a_{\ell}r_2^{-\ell}\sum_{k=j-\ell}^{j-1}\left(\frac{r_1}{r_2}\right)^k {k+\ell\choose M_2} k^{M_1}\nonumber\\
        &= \ \sum_{\ell=1}^{j}a_{\ell}r_2^{-j}r_1^{-\ell}\sum_{k=0}^{\ell-1}\left(\frac{r_1}{r_2}\right)^k{j+k\choose M_2}(j+k-\ell)^{M_1},
    \end{align}
    we see that it suffices to show that $b_{\ell}\defeq r_2^{-j}r_1^{-\ell}\sum_{k=0}^{\ell-1}\left(r_1/r_2\right)^k{j+k\choose M_2}(j+k-\ell)^{M_1}$ satisfies the linear recurrence with characteristic polynomial $z^j p(1/z)=a_jz^j+a_{j-1}z^{j-1}+\cdots+a_1z-1$ for all $\ell\geq j$. We proceed by determining the generating function of $\{b_{\ell}\}_{\ell\geq0}$. Splitting into two series, we observe that
    \begin{align}\label{b_l_comp_1}
        \sum_{\ell\geq0}b_{\ell}z^\ell& \ = \ \sum_{\ell\geq0}z^{\ell}\left(r_2^{-j}r_1^{-\ell}\sum_{k=0}^{\ell}\left(\frac{r_1}{r_2}\right)^k{j+k\choose M_2}(j+k-\ell)^{M_1}\right)-\sum_{\ell\geq0}z^{\ell}\left(r_2^{-j}r_1^{-\ell}\left(\frac{r_1}{r_2}\right)^{\ell}{j+\ell\choose M_2}j^{M_1}\right)\nonumber\\
        & \ = \ r_2^{-j}\sum_{\ell\geq0}\left(\frac{z}{r_1}\right)^{\ell}\sum_{k=0}^{\ell}\left(\frac{r_1}{r_2}\right)^k{j+k\choose M_2}(j+k-\ell)^{M_1}-\frac{D(z)}{(1-z/r_2)^{M_2+1}}
    \end{align}
    for some polynomial $D(z)$ of degree at most $M_2$.
    Recognizing that $\sum_{k=0}^{\ell}\left(r_1/r_2\right)^k{j+k\choose M_2}(j+k-\ell)^{M_1}$ is a Cauchy product, we see that
    \begin{align}\label{b_1_comp_2}
        r_2^{-j}\sum_{\ell\geq0}\left(\frac{z}{r_1}\right)^{\ell}\sum_{k=0}^{\ell}\left(\frac{r_1}{r_2}\right)^k{j+k\choose M_2}(j+k-\ell)^{M_1} \ = \ \frac{A(z)}{(1-z/r_1)^{M_1+1}(1-z/r_2)^{M_2+1}}
    \end{align}
    for some polynomial $A(z)$ of degree at most $M_1+M_2$. Let $L(z)=(1-z/r_1)^{M_1+1}(1-z/r_2)^{M_2+1}$. Since $r_1$ and $r_2$ are distinct, $\deg(A)\leq M_1+M_2<j$. So \eqref{b_l_comp_1} and \eqref{b_1_comp_2} show that $\{b_{\ell}\}_{\ell\geq0}$ is a linear combination of sequences which, by \cite[pg.255]{analytic_combo}, all satisfy the linear recurrence with characteristic polynomial $z^{M_1+M_2+2}L(1/z)$ for all $\ell\geq j$. As $r_1,r_2$ are distinct, $z^{M_1+M_2+2} L(1/z)$ divides $z^j p(1/z)$, whence it follows that $b_{\ell}$ satisfies the linear recurrence with characteristic polynomial $z^j p(1/z)$ for all $\ell\geq j$.
\end{proof}

\begin{definition}
    Let $\mathbf{u}=\{c_n\}$ be as in Theorem \ref{main_theorem}. For each characteristic root $r$, let $p^{\mathbf{u}}_r(n)$ be the polynomial such that $c_n=\sum_{r}p^{\mathbf{u}}_{r}(n) r^n$ for all $n\geq0$.
\end{definition}

\begin{proof}[\textbf{Proof of Theorem \ref{main_theorem}}]
    By a standard result in the theory of recurrence relations (e.g., \cite[Theorem IV.9]{analytic_combo}), we know that $c_n=\sum_r\sum_{i=1}^{m_r}\kappa_{r,i}n^{m_r-i}r^n$
    for some constants $\kappa_{r,i}$. So it suffices to show that formula \eqref{subset_formula_kappa} holds for all characteristic roots $r$ and $1\leq i\leq m_r$.
    
    Let $\mathbf{u}=\{c_n\}_{n\geq0}$. Let $r$ be a characteristic root, let $1\leq i\leq m_r$, and let $S\not\ni r$ be the set of characteristic roots greater than or equal to $r$ in modulus. Define
    \begin{align}
        \mathbf{u}' \ = \ \left\{c_n-\sum_{R\in S}p^{\mathbf{u}}_R(n) R^n-\sum_{v=1}^{i-1}\kappa_{r,v} n^{m_r-v} r^n\right\}_{n\geq0}.
    \end{align}
    Noting that $\kappa_{r,i}=[z^{m_r-i}]~p_r^{\mathbf{u}}(z)=[z^{m_r-i}]~p_r^{\mathbf{u'}}(z)$, we proceed by computing $[z^{m_r-i}]~p_r^{\mathbf{u'}}(z)$. As $n^M R^n$ satisfies recurrence \eqref{c_n_recurrence} for any characteristic root $R$ and $M\leq m_R-1$, it follows that $\mathbf{u}'$ also satisfies recurrence \eqref{c_n_recurrence} with initial conditions $\{c_k-\sum_{R\in S}p^{\mathbf{u}}_R(k) R^k-\sum_{v=1}^{i-1}\kappa_{r,v} k^{m_r-v} r^k\}_{0\leq k\leq j-1}$. Let $p(z)$ be the characteristic polynomial of \eqref{c_n_recurrence}. Then by Proposition \ref{generating_function_formula}, the generating function of $\mathbf{u}'$ is
    \begin{align}
        f(z) \ \defeq \ \frac{1}{p(z)}\left(\sum_{k=0}^{j-1}\left(c_k-\sum_{R\in S}p^{\mathbf{u}}_R(k) R^k-\sum_{v=1}^{i-1}\kappa_{r,v} k^{m_r-v} r^k\right)\sum_{\ell\geq j-k}a_{\ell}z^{k+\ell}+p(z)\sum_{k=0}^{j-1}c_k z^k\right).
    \end{align}
    Note that $f$ has a single dominant pole $r^{-1}$ with multiplicity $m_r-i+1$ since the asymptotically dominant term of $\mathbf{u}'$ is $\kappa_{r,i}n^{m_r-i}r^n$. So by \cite[Note IV.26]{analytic_combo}, $\kappa_{r,i}=C/(m_r-i)!$, where $C=\lim_{z\to r^{-1}}(1-rz)^{m_r-i+1}f(z)$. Since the numerator of $f$ is divisible by $(1-rz)^{i-1}$, it follows that $C=C_1/C_2$, where $C_2=\lim_{z\to r^{-1}}p(z)/(1-rz)^{m_r}$ and $C_1$ equals
    \begin{align}\label{def_of_C_1}
        \lim_{z\to r^{-1}}\frac{1}{(1-rz)^{i-1}}\left(\sum_{k=0}^{j-1}\left(c_k-\sum_{R\in S}p^{\mathbf{u}}_R(k) R^k-\sum_{v=1}^{i-1}\kappa_{r,v} k^{m_r-v} r^k\right)\sum_{\ell\geq j-k}a_{\ell}z^{k+\ell}+p(z)\sum_{k=0}^{j-1}c_k z^k\right).
    \end{align}
    As $i-1< m_r$, 
    \begin{align}
        \lim_{z\to r^{-1}}\frac{p(z)}{(1-rz)^{i-1}}\sum_{k=0}^{j-1}c_k z^k \ = \ 0;
    \end{align}
    hence, \eqref{def_of_C_1} reduces to
    \begin{align}
        \lim_{z\to r^{-1}}\frac{1}{(1-rz)^{i-1}}\left(\sum_{k=0}^{j-1}\left(c_k-\sum_{R\in S}p^{\mathbf{u}}_R(k) R^k-\sum_{v=1}^{i-1}\kappa_{r,v} k^{m_r-v} r^k\right)\sum_{\ell\geq j-k}a_{\ell}z^{k+\ell}\right).
    \end{align}
    Thus,
    \begin{align}
        C_1 \ &= \ \frac{(-1)^{i-1}}{r^{i-1}(i-1)!}\sum_{k=0}^{j-1}\left(c_k-\sum_{R\in S}p^{\mathbf{u}}_R(k) R^k-\sum_{v=1}^{i-1}\kappa_{r,v} k^{m_r-v} r^k\right)\sum_{\ell\geq j-k}(i-1)!{k+\ell\choose i-1}a_{\ell}r^{-k-\ell+i-1}\nonumber\\
        &= \ (-1)^{i-1}\sum_{k=0}^{j-1}\left(c_k-\sum_{R\in S}p^{\mathbf{u}}_R(k) R^k-\sum_{v=1}^{i-1}\kappa_{r,v} k^{m_r-v} r^k\right)\sum_{\ell\geq j-k}{k+\ell\choose i-1}a_{\ell}r^{-k-\ell}\nonumber\\
        &= \ (-1)^{i-1}\sum_{k=0}^{j-1}\left(c_k-\sum_{R\in S}p^{\mathbf{u}}_R(k) R^k-\sum_{v=1}^{i-1}\kappa_{r,v} k^{m_r-v} r^k\right)s_{r,k,i-1}\nonumber\\
        &\overset{*}{=} \ (-1)^{i-1}\sum_{k=0}^{j-1}\left(c_k-\sum_{v=1}^{i-1}\kappa_{r,v} k^{m_r-v} r^k\right)s_{r,k,i-1},
    \end{align}
    where the starred equality follows from Lemma \ref{zeroing_lemma}. Computing $C_2$, we obtain
    \begin{align}
        C_2 \ = \  \frac{-\sum_{\ell\geq m_r}m_r!{\ell\choose m_r}a_{\ell}r^{-\ell+m_r}}{(-1)^{m_r}r^{m_r}m_r!} \ = \ (-1)^{m_r-1}\sum_{\ell\geq m_r}{\ell\choose m_r}a_{\ell}r^{-\ell}=N_r^{-1}.
    \end{align}
    Therefore, we have the recurrence
    \begin{align}\label{kappa_recurrence}
        \kappa_{r,i} \ = \ \frac{(-1)^{i-1}N_r}{(m_r-i)!}\sum_{k=0}^{j-1}\left(c_k-\sum_{v=1}^{i-1}\kappa_{r,v}k^{m_r-v}r^k\right)s_{r,k,i-1}.
    \end{align}
    It now only remains to show that recurrence \eqref{kappa_recurrence} admits the explicit formula
    \begin{align}\label{kappa_formula}
        \kappa_{r,i} \ = \ \frac{(-1)^{i-1}}{(m_r-i)!}\sum_{\substack{\varnothing\neq\{d_1,\dots,d_p\}\subseteq\{0,1,\dots,i-1\}\\  d_1<\cdots<d_p=i-1}}{N}_r^p\sum_{k_1=0}^{j-1}c_{k_1}s_{r,k_1,d_1}\prod_{u=1}^{p-1}\sum_{k_2=0}^{j-1}\frac{k_2^{m_r-1-d_u}r^{k_2} s_{r,k_2,d_{u+1}}}{(-1)^{d_u-1}(m_r-1-d_u)!}.
    \end{align}
    Let $\{\kappa_{r,i}'\}_{i\in\ZZ^+}$ be the sequence defined by the explicit formula given in \eqref{kappa_formula}. For abbreviation, let 
    \begin{align}
        \Pi(\{d_1,\dots,d_p\}) \ = \ \sum_{k_1=0}^{j-1}c_{k_1}s_{r,k_1,d_1}\prod_{u=1}^{p-1}\sum_{k_2=0}^{j-1}\frac{k_2^{m_r-1-d_u}r^{k_2} s_{r,k_2,d_{u+1}}}{(-1)^{d_u-1}(m_r-1-d_u)!}.
    \end{align}
    Towards proving that $\kappa_{r,i}\equiv\kappa_{r,i}'$, it suffices to show that $\{\kappa_{r,i}'\}_{i\in\ZZ^+}$ satisfies the recurrence since $\{\kappa_{r,i}\}_{i\in\ZZ^+}$ is completely determined by \eqref{kappa_recurrence}. Note
    \begin{align}
        \sum_{v=1}^{i-1}\kappa_{r,v}'\sum_{k=0}^{j-1}{k^{m_r-v}r^k s_{r,k,i-1}} \ = \ &\sum_{v=1}^{i-1}\sum_{\substack{\varnothing\neq\{d_1,\dots,d_p\}\subseteq\{0,1,\dots,v-1\}\\  d_1<\cdots<d_p=v-1}}{N}_r^{p}~\Pi(\{d_1,\dots,d_p\})\sum_{k=0}^{j-1}\frac{k^{m_r-v}r^k s_{r,k,i-1}}{(-1)^{v-1}(m_r-v)!}\nonumber\\
        = \ &-\sum_{v=1}^{i-1}\sum_{\substack{\varnothing\neq\{d_1,\dots,d_{p+1}\}\subseteq\{0,1,\dots,i-1\}\\  d_1<\cdots<d_p=v-1,~ d_{p+1}=i-1}}{N}_r^{p}~\Pi(\{d_1,\dots,d_{p+1}\}).\label{kappa_prime_sum}
    \end{align}
    Now with the aid of \eqref{kappa_prime_sum}, we have
    \begin{align}
        &\frac{(-1)^{i-1}N_r}{(m_r-i)!}\sum_{k=0}^{j-1}\left(c_k-\sum_{v=1}^{i-1}\kappa'_{r,v}k^{m_r-v}r^k\right)s_{r,k,i-1}\nonumber\\
        = \ &\frac{(-1)^{i-1}}{(m_r-i)!}\Bigg(N_r\sum_{k=0}^{j-1}c_ks_{r,k,i-1}
        -N_r\sum_{v=1}^{i-1}\kappa_{r,v}'\sum_{k=0}^{j-1}k^{m_r-v}r^ks_{r,k,i-1}\Bigg)\nonumber\\
        = \ &\frac{(-1)^{i-1}}{(m_r-i)!}\left(N_r\Pi(\{i-1\})+\sum_{v=1}^{i-1}\sum_{\substack{\varnothing\neq\{d_1,\dots,d_{p+1}\}\subseteq\{0,1,\dots,i-1\}\\  d_1<\cdots<d_p=v-1,~ d_{p+1}=i-1}}{N}_r^{p+1}~\Pi(\{d_1,\dots,d_{p+1}\})\right)\nonumber\\
        = \ & \frac{(-1)^{i-1}}{(m_r-i)!}\sum_{\substack{\varnothing\neq\{d_1,\dots,d_{p+1}\}\subseteq\{0,1,\dots,i-1\}\\  d_1<\cdots<d_{p+1}=i-1}}{N}_r^{p+1}~\Pi(\{d_1,\dots,d_{p+1}\}) \ = \ \kappa'_{r,i}.
    \end{align}
    Thus, $\kappa_{r,i}\equiv\kappa_{r,i}'$.
\end{proof}

\section{Directions for Future Work}

A natural continuation of our work is to generalize Theorem \ref{main_theorem} to a broader class of recurrences. A particularly natural generalization is to consider linear recurrences of the form $c_n=a_1 c_{n-1}+\cdots+a_j a_{n-j}+d(n)$, where $d(n)$ is some polynomial. These recurrences are classified as non-homogeneous linear recurrences: a linear recurrence $c_n-a_1 c_{n-1}-\cdots-a_j a_{n-j}=f(n)$ is \emph{homogeneous} if $f(n)\equiv0$, and \emph{non-homogeneous} otherwise. This proposed generalization only encapsulates the case when $f$ is a polynomial; however, this induces the generating function of $\{c_n\}$ to remain rational, thereby permitting the techniques presented in this work to be readily generalized.

\section*{Acknowledgments}

This work was supported by the National Science Foundation
grant \textit{DMS-2241623}, Williams College, the Finnerty Fund, and the Winston Churchill Foundation. We are grateful to the anonymous referee for their feedback, which prompted a vast improvement of a previous version of this work. We thank Jay Pantone, Robin Pemantle, Richard Stanley, Vincent Vatter, and Doron Zeilberger for their consultation on this project. We acknowledge Joe Cooper for helpful conversations and for highlighting errors and typos in an earlier draft. We also thank participants of the 21\textsuperscript{st} International Fibonacci Conference for comments on an earlier version of this work.
    
\bibliography{REFS_recurrence_relations_SMALL_2024_04.tex}{}
\bibliographystyle{plain}

\end{document}